\documentclass[12pt]{article}
\usepackage[dvips]{graphics}
\usepackage{amsmath,amssymb}
\begin{document}

\newcommand{\de}{\mbox{$\delta $}}
\newcommand{\e}{\mbox{${\bf \epsilon \rm}$}}
\newcommand{\qed}{\mbox{$ \quad  \qquad \qquad \qquad \qquad \qquad \qquad \qquad \qquad \qquad \qquad
 \qquad \qquad \qquad \qquad \quad \qquad \square $ }}
\newtheorem{cor}{Corollary}[section]
\newtheorem{alg}{Algorithm}[section]
\newtheorem{lemma}{Lemma}[section]
\newtheorem{theo}{Theorem}[section]
\newtheorem{defi}{Definition}[section]
\newtheorem{exa}{Example} [section]
\newtheorem{pro}{Proposition}[section]
\newtheorem{rem}{Remark}[section]
\noindent {\bf AN ALGORITHM FOR DESIGNING FEEDBACK STABILIZERS
OF NONLINEAR POLYNOMIAL SYSTEMS. \rm}
\vskip 10 pt \noindent {\it Stelios Kotsios }
\par
\noindent {\it University of Athens, Faculty of Economics, Department of Mathematics
and Computer Science. Pesmazoglou 8, 10559 Athens, Greece } \vskip 15 pt \noindent
{\bf Abstract:}{\it The aim of this paper is to present a symbolic computational
algorithm that will allow us to deal with the feedback stabilization problem for
continuous nonlinear polynomial systems.
 The overall approach is based on a methodology that checks 
 the positivity of a given polynomial.}

\section{Introduction}
One of the major purposes of control theory is the analysis and the design of feedback
control systems. A great number of both analytical and computational methods have also
been adopted for continuous and discrete systems, \cite{kn:isidori},\cite{kn:sontag}.
\par
\noindent In this paper we examine the feedback stabilization problem for a large
variety of continuous nonlinear systems. Our methodology resembles that of the applied
computational tools used to solve control problems \cite{kn:karkanias}. Specifically,
we consider a system of the form $\dot{\bf x}={\bf \Phi}({\bf x}, {\bf u})$,
$(\Sigma)$, where the components of ${\bf \Phi}$, $\Phi_i$, are multivariable
polynomials. According to the theory, if the linearization of $(\Sigma)$ at the origin
is asymptotically controllable (i.e., all the uncontrollable eigenvalues have negative
real parts), then the nonlinear system is locally asymptotically stabilizable
\cite{kn:sontag}. This means that we can find a feedback-law that makes the
closed-loop system asymptotically stable to the origin. The aim of this paper is to
calculate those feedback-laws computationally by means of certain symbolic algorithms.
Concretely, we seek for nonlinear state-feedbacks of the form ${\bf
u}={\bf a}({\bf x})$, where ${\bf a}({\bf x})$ consists of multivariable polynomials,
so that a proper Lyapunov function becomes negative. To solve the problem we develop
the following algorithms:
\par
\noindent {\bf The Formal Algorithm}. This algorithm allows us to write a polynomial
$p$ as follows:
\[ p=c_{1}(W_{i,\sigma,\varphi}) x_1^{j_{1,1}} \cdot
[ W_{2,1,1}x_1+x_2]^{j_{2,1}} \cdot [W_{3,1,1} x_1+W_{3,2,1} x_2 +x_3]^{j_{3,1}}
\cdots \]
\[\cdots [W_{n,1,1}x_1 +W_{n,2,1}x_2+ \cdots +x_n]^{j_{n,1}}+ \]
\[ +c_{2}(W_{i,\sigma,\varphi}) x_1^{j_{1,2}} \cdot
[ W_{2,1,2}x_1+x_2]^{j_{2,2}} \cdot [W_{3,1,2} x_1+W_{3,2,2} x_2 +x_3]^{j_{3,2}}
\cdots \]
\[\cdots [W_{n,1,2}x_1 +W_{n,2,2}x_2+ \cdots +x_n]^{j_{n,2}}+ \cdots \]
\[+c_{k}(W_{i,\sigma,\varphi}) x_1^{j_{1,k}} \cdot
[ W_{2,1,k}x_1+x_2]^{j_{2,k}} \cdot [W_{3,1,k} x_1+W_{3,2,k} x_2 +x_3]^{j_{3,k}}
\cdots \]
\begin{equation}\label{intro}
\cdots [W_{n,1,k}x_1 +W_{n,2,k}x_2+ \cdots +x_n]^{j_{n,k}}+ R_{\cal W}(x_1)
\end{equation}
where the exponents $j_{a,b}$ are specific positive whole numbers, the quantities
$W_{i,\sigma,\varphi}$ are undetermined parameters that can take certain values,
$c_j(W_{i,\sigma,\phi})$ the coefficients depending on the parameters
$W_{i,\sigma,\phi}$ and $R_{\cal W}(x_1)$ a polynomial of the single variable $x_1$,
called the remainder. Equation (\ref{intro}) is called the
Formal-Linear-Like-Factorization of $p$ and appeared firstly in
\cite{kn:hercma}. The essential tool of this methodology is a
continuous reduction of $p$, by means of an Euclidean division. This makes our method
similar to others \cite{kn:cox}, \cite{kn:parrilo}. Yet in our methodology we only
deal with a specific polynomial, in addition to the Gr\"{o}bner basis that works with a
polynomial ideal.
\par
\noindent {\bf The POS-Algorithm.} This algorithm checks the positivity of a
multivariable polynomial $p$. This task can be achieved by giving to
$W_{i,\sigma,\phi}$ values such that all the non-square terms to be eliminated and the
square terms have positive coefficients. A lot of work has been done in this direction
and some of the results can be found in the following
\cite{kn:lasserre},\cite{kn:hasan},\cite{kn:prestel}, to mention but a few. In recent
years there has been a strong algorithmic approach to the whole issue, as can be seen,
for instance, in \cite{kn:hochbaum},\cite{kn:mishra},\cite{kn:boyd}. A number of those
algorithms dealt with the problem of the existence of a sum of squares decomposition
\cite{kn:choi},\cite{kn:powers},\cite{kn:parriloPHD}. That is when we are able to
write a polynomial $p$ as $p=\sum_{i}f_i^2(x_1,x_2,\ldots,x_n)$. In our sum of squares,
$p=\sum_{i}f_i^2(x_1,x_2,\ldots,x_n)$, the quantities $f_i(x_1,x_2,\ldots,x_n)$ are "
linear " polynomials of the variables $x_1,x_2,\ldots,x_n$, while in other approaches
\cite{kn:parriloPHD}, $f_i(x_1,x_2,\ldots,x_n)$ can be arbitrary. Moreover, our
methodology provides not a specific sum of squares but a whole class, depending on the
various proper values of the parameters $W_{i,\sigma,\varphi}$.
\par
\noindent {\bf The Feedback-GAS-Algorithm}. This algorithm accepts as input the
polynomials $\Phi_i$ and a Lyapunov Function $L$, and calculates the feedback
connection ${\bf u}={\bf a}({\bf x })$. To achieve this we assume that ${\bf a}({\bf
x})$ consist of a multivariable polynomial with parametric coefficients (we denote
them by $A^{(j)}_{i_1,i_2,\ldots,i_n}$). Our aim is to determine those values of the
parameters $W_{i,\sigma,\varphi}$ and $A^{(j)}_{i_1,i_2,\ldots,i_n}$ that make the
Lyapunov function negative along the trajectories of the closed-loop system. This
guarantees the stability of the origin. The main merits of our method are the
following:
\par
\noindent 1) It constitutes a pragmatic computational method. Indeed, the
feedback-laws are derived from symbolic computational algorithms. Appropriate software
has been created for this purpose, and all the examples presented in the current paper
have been studied with the aid of this software.
\par
\noindent 2) It provides us with a whole class of admissible controllers, as opposed
to only a single one.
\par
\noindent 3) It works even when the linearization of the nonlinear system has an
uncontrollable eigenvalue with zero real part, \cite{kn:sontag}.
\par
\noindent 4) Since the Lyapunov function, upon construction, is negative for every
point, the stability is global. \vskip 8 pt \noindent The paper is divided into two
parts. The first part is devoted to the presentation of the algorithms used throughout
the paper. The second part is devoted to the specific control application. ${\bf R}$
and ${\bf Z}^+$ will denote the sets of real numbers and positive integers
correspondingly.

\section{The Algebraic Background}
Let ${\bf R}$ be the set of real numbers and $x_1,x_2,\ldots,x_n$ $n$-variables. An
expression of the form $ p=\sum_{\lambda=1}^{\varphi} c_{\lambda} x_1^{a_{1,\lambda}}
x_2^{a_{2,\lambda}} \cdots  x_n^{a_{n,\lambda}}$, where $c_{\lambda} \in {\bf R}$ and
some of the exponents $a_{i,j} \in \mathbf{Z}^+$ are not equal to zero, is called a
polynomial in $x_1,x_2,\ldots,x_n$ with real coefficients or, for short, a real
polynomial. The set of all real polynomials in $x_1,x_2,\ldots,x_n$ is denoted by
${\bf R}[x_1,x_2,\ldots,x_n]$. An element $x_1^{a_{1,\lambda}} x_2^{a_{2,\lambda}}
\cdots  x_n^{a_{n,\lambda}}$ is called a {\it monomial} and an element $c_{\lambda}
x_1^{a_{1,\lambda}} x_2^{a_{2,\lambda}} \cdots  x_n^{a_{n,\lambda}}$ is called a {\it
term}. Let $\phi_{n,\lambda}= x_1^{a_{1,\lambda}} x_2^{a_{2,\lambda}} \cdots
x_n^{a_{n,\lambda}}$ and $\phi_{m,\mu}= x_1^{a_{1,\mu}} x_2^{a_{2,\mu}} \cdots
x_m^{a_{m,\mu}}$ be two monomials. We defined the {\it lexicographical } order among
monomials \cite{kn:cox}, as follows: we say that $\phi_{n,\lambda}$ is ordered less
than $\phi_{m,\mu}$, denoted by $\phi_{n,\lambda} \prec \phi_{m,\mu}$, if either $n<
m$ or $n=m$ and $a_{n,\lambda}<a_{m,\mu}$. In other words, the monomials are ordered
as follows: $ x_1\prec \cdots \prec x_1^7 \prec$ $ \cdots \prec x_1x_2 \prec$ $ \cdots
\prec x_1 x_2^8 \prec$ $ \cdots \prec x_1x_2x_3 \prec \cdots$ Let $p$ be a given
polynomial. Ordered lexicographically, the term that corresponds to the maximum
monomial is called the {\it maximum term} denoted by $maxterm(p)$. Throughout the
paper, a variable $W_{i,j,k}$ taking values in ${\bf R}$ is called {\it an undetermined
parameter}. The set of undetermined parameters is denoted by ${\cal W}$.
\par
\noindent Let $p$ $\in \mathbf{R}$$ [{ \rm x_1,x_2,\ldots,x_n} ]$ be a polynomial with
$n$-variables and ${\cal W}=\{ W_{i,\sigma,\varphi} \}$ a set of undetermined
parameters, taking values in $\mathbf{R}$. A Formal-Linear-Like-Factorization of $p$
is an expression of the form:
\[ p=\sum_{\mu =1}^k c_{\mu}(W_{i,\sigma,\varphi}) x_1^{j_{1,\mu}} \cdot
[ W_{2,1,\mu}x_1+x_2]^{j_{2,\mu}} \cdot [W_{3,1,\mu} x_1+W_{3,2,\mu} x_2
+x_3]^{j_{3,\mu}} \cdots \]
\[ \cdots [W_{n,1,\mu}x_1 +W_{n,2,\mu}x_2+ \cdots +x_n]^{j_{n,\mu}}+R_{\cal W}(x_1) \]
where the coefficients $c_{\mu}(W_{i,\sigma,\varphi})$ are polynomial functions of the
parameters $W_{i,\sigma,\varphi}$ and the {\it remainder} $R_{\cal W}(x_1)$ is a
polynomial only of the single variable $x_1$, with coefficients depending on the
parameters $W_{i,\sigma,\varphi}$, too. Some of the exponents $j_{1,\mu},j_{2,\mu},
\ldots , j_{n,\mu}$ $\in {\bf Z^+}$, may be equal to zero. 
The {\it Formal-Linear-Like-Factorization} of $p$
is denoted by {\it FormalLF[p]}.  In a number of instances,  a
Formal-Linear-Like-Factorization of $p$, is written for short as $
FormalLF[p]=\sum_{\mu =1}^k c_{\mu}L^{j_{1,\mu}}_{1,\mu} \cdot
 L^{j_{2,\mu}}_{2,\mu} \cdot L^{j_{3,\mu}}_{3,\mu} \cdots L^{j_{n,\mu}}_{n,\mu}+R$,
where $ L_{\sigma,\mu}=W_{\sigma,1,\mu}x_1+W_{\sigma,2,\mu}x_2+ \cdots + x_{\sigma}$,
$\sigma=1, \ldots, n$ and $c_{\mu}$,$R$ are the abbreviations of
$c_{\mu}(W_{i,\sigma,\varphi})$ and $R_{\cal W}(x_1)$.

\begin{exa}
We have the polynomial $p=5x_1-7x_1x_2+11x_1x_3$. The Formal-Linear-Like-Factorization
of $p$ is, $ p=11x_1(W_{3,1,1}x_1+$ $W_{3,1,2}x_2+x_3)+$ $ (-7-11W_{3,1,2})$
$x_1(W_{2,1,2}x_1+x_2)+$ $5x_1+(7W_{2,1,2}-$ $11W_{3,1,1}+$
$11W_{2,1,2}W_{3,2,1})x_1^2$.
 \end{exa}
\noindent The following theorem deals with the question of uniqueness of the
Formal-Linear-Like-Factorization.

\begin{theo}
For a given polynomial $p \in  \mathbf{R} [{\rm x_1,x_2,\ldots,x_n} ]$, the
$FormalLF[p]$ is unique, under the assumption that the parameters ${\cal W}=\{
W_{i,j,\varphi} \}$ are considered as constants.

\end{theo}
{\bf Proof:} If $p$ has only $x_1$-terms, the proof of the theorem is trivial, with
$c_{\mu}=0$ and $R=p$. Let us suppose that $p$ has at least one term other than the
$x_1$-terms. Let us further suppose that $\lambda x_1^{j_{1,h}}x_2^{j_{2,h}} \cdots
x_n^{j_{n,h}}$ is its maximum term. This term appears also in the product:
$c_hL_{1,h}^{j_{1,h}}L_{2,h}^{j_{2,h}} \cdots L_{n,h}^{j_{n,h}}$. By equating their
coefficients, we calculate the quantity $c_h$, uniquely, actually $c_h=\lambda$.
Repeating the same procedure for the term with the next higher order, we find an
expression for the "next" coefficient $c_{h-1}$. Since this expression is a function
of $c_h$ and some of the parameters $W_{i,\sigma,\varphi}$ are considered as
constants, we conclude that $c_{h-1}$ is also defined uniquely. By induction, we
finally get coefficients $c_{\mu}, \mu=1, \ldots,k$, all of which are uniquely
determined. The polynomial $R$ consists only of $x_1$-terms. These terms arise either
from the polynomial $p$ or from the products ${c}_{\mu}L^{j_{1,\mu}}_{1,\mu} \cdot
 L^{j_{2,\mu}}_{2,\mu} \cdot L^{j_{3,\mu}}_{3,\mu} \cdots L^{j_{n,\mu}}_{n,\mu}$. The unique determination of the coefficients
 ${c}_{\mu}$ entails the uniqueness of $R$, and the theorem has been proved.
$\square $
\par
\noindent What is of interest is the issue of constructing the
Formal-Linear-Like-Factorization of a given polynomial. This we can do through the
algorithm we describe below. It is based on a methodology that is analogous to the
usual Euclidean division among polynomials and can be implemented on a computer via
proper software (MATHEMATICA, for instance).
\par
\vskip 10 pt \noindent \underline{\bf THE FORMAL ALGORITHM} \vskip 10 pt \small {\sf
\par
\noindent {\bf Input:} A multivariable polynomial $p$, a set of undetermined parameters
${\cal W}=\{ W_{i,\sigma,\varphi}\}$, taking values in ${\bf R}$.
\par
\noindent {\bf Initial Conditions:} $k=0$
\par
\begin{description}
\item{\bf Step 1:} We set $k=k+1$.

\item{\bf Step 2:} We find the maximum term of $p$,
$maxterm(p)=c_k x_1^{j_{1,k}}\cdots x_n^{j_{n,k}}$. The coefficient $c_k$, in the
first iteration, is a constant number. Then, it depends on the set of parameters
${\cal W}$.

\item{\bf Step 3:} We form the linear polynomials:
$ L_{1,k}= x_1$, $ L_{2,k}=W_{2,1,k}x_1+x_2$, $\ldots$,
$L_{n,k}=W_{n,1,k}x_1+W_{n,2,k}x_2+W_{n,3,k}x_3+\cdots +x_n$.

\item{\bf Step 4:} We make the subtraction:
$ R_k=p-c_{k}L_{1,k}^{j_{1,k}} L_{2,k}^{j_{2,k}}\cdots L_{n,k}^{j_{n,k}}$.

\item{\bf Step 5:} {\bf IF } $R_k$ depends only on the variable $x_1$
{\bf THEN} set $R=R_k$ and go to the output {\bf ELSE} put $p=R_k$ and go to step 1.
\end{description}
 {\bf Output:}
The quantities $R$, $c_\mu$, $L_{i,\mu}$, $j_{i,\mu}$, $\mu=1, \ldots, k$, $i=1,
\ldots, n$ }

\vskip 15 pt \normalsize \noindent The following theorem proves the finiteness and the
efficiency of the algorithm.
\begin{theo}
The Formal Algorithm terminates after a finite number of steps. If $R, c_\mu,
L_{i,\mu}$, $j_{i,\mu}$, $\mu=1, \ldots,k$, $i=1, \ldots, n$ are its outputs, then:
\[ FormalLF[p]=\sum_{\mu =1}^k c_{\mu}L^{j_{1,\mu}}_{1,\mu} \cdot
 L^{j_{2,\mu}}_{2,\mu} \cdot L^{j_{3,\mu}}_{3,\mu} \cdots L^{j_{n,\mu}}_{n,\mu}+R\]
 \end{theo}
 {\bf Proof:}
Let $p$ be a multivariable polynomial and $z=\gamma x_1^{a_1}x_2^{a_2} \cdots
x_n^{a_n}$ its maximum term. We follow the Formal Algorithm step by step. When $k=1$,
step 2 will give $c_1=\gamma$ and $j_{1,1}=a_1, \ldots, j_{n,1}=a_n$. Taking into
consideration the above values and the construction of the linear polynomials
$L_{1,1},L_{2,1}, \ldots, L_{n,1}$, step 4 will produce a new polynomial $R_1$, which
will not contain the term $z$. Obviously, the maximum term of $R_1$, which is called
$z_1$, will be ordered lower than $z$, $z_1 \prec z$, with respect to the order raised
earlier. By induction, we get for the maximum terms $z\succ z_1 \succ z_2 \succ z_3
\succ \cdots $. This nest and the construction of the order will finally eliminate all
but $x_1$-terms. This fact guarantees the termination of the algorithm. Substituting
now reversely, we have successively $ p=c_1 L_{1,1}^{j_{1,1}} \cdots L_{n,1}^{j_{n,1}}
+R_1$, $R_1=c_2 L_{1,2}^{j_{1,2}} \cdots L_{n,2}^{j_{n,2}} +R_2$, $ \cdots $, $
R_{k-1}=c_k L_{1,k}^{j_{1,k}} \cdots L_{n,k}^{j_{n,k}} +R$. Combining these we get: $
p= \sum_{\mu =1}^k c_\mu L_{1,\mu}^{j_{1,\mu}} \cdots L_{n,\mu}^{j_{n,\mu}} +R$, which
is the Formal-Linear-Like-Factorization upon request. $ \square $ 
\par \noindent The following example
exhibits the function of the algorithm.
\begin{exa}\label{ex1}
We have the previous polynomial $p=5x_1-7x_1x_2+11x_1x_3$. We want to find a
Formal-Linear-Like-Factorization of $p$. In order to clarify our ideas we shall follow
the Formal-Algorithm in detail. First, $maxterm(p)=11x_1x_3$, here $c_1=11$ and
$j_{1,1}=1, j_{1,2}=0, j_{1,3}=1$. The linear polynomials $L_{i,1}$ are $L_{1,1}=x_1$,
$L_{2,1}=W_{2,1,1}x_1+x_2$, $L_{3,1}=W_{3,1,1}x_1+W_{3,1,2}x_2+x_3$. Then we get $
R_1=p(x_1,x_2,x_3)-11 L_{1,1}^1 \cdot L_{2,1}^0 \cdot L_{3,1}^1=$
$5x_1-7x_1x_2-11x_1(W_{3,1,1}x_1+W_{3,1,2}x_2+x_3)=$
$5x_1-11W_{3,1,1}x_1^2+(-7-11W_{3,1,2})x_1x_2$. The new maxterm is
$(-7-11W_{3,1,2})x_1x_2$, with $c_2=-7-11W_{3,1,2}$, $j_{1,2}=1$, $j_{2,2}=1$. The
polynomials $L_{i,2}$ now become $L_{1,2}=x_1$, $L_{2,2}=W_{2,1,2}x_1+x_2$ and $R_2$
 is,
$R_2=R_1-(-7-11W_{3,1,2})x_1(W_{2,1,2}x_1+x_2)=$
$5x_1+(7W_{2,1,2}-11W_{3,1,1}+11W_{2,1,2} W_{3,2,1})x_1^2$.
 $R_2$ contains only $x_1$-terms and the algorithm terminates. Therefore, the Formal
 Linear-Like Factorization of $p$ is
$ p=11x_1(W_{3,1,1}x_1+W_{3,1,2}x_2+x_3)+$ $ (-7-11W_{3,1,2})x_1(W_{2,1,2}x_1+x_2)+$
$+ 5x_1+(7W_{2,1,2}-11W_{3,1,1}+11W_{2,1,2} W_{3,2,1})x_1^2$
\end{exa}
\noindent We can take different expressions of a concrete polynomial $p$, by giving
certain values to the parameters $W_{i,\sigma,\varphi}$. Such procedures are called
{\it evaluations} of the $FormalLF[p]$. The most rigorous approach is the following:
Let ${\cal W}=\{W_{i,\sigma, \varphi}\}$ be the set of the variables that appear in the
Formal - Linear - Like - Factorization of a given polynomial $p$. By arranging the
parameters in an increasing order we form the vector ${\cal
W}=(W_{i_k,\sigma_k,\varphi_k})_{k=1,2, \ldots,n}$. Let ${\bf
r}=(a_k)_{k=1,2,\ldots,n}$ be a vector of real numbers which has the same length as
the vector ${\cal W}$. We say that the parameters ${\cal W}$ follow the rules ${\bf
r}$ and we write ${\cal W} \to  {\bf r}$ if the following substitution is valid:
$W_{i_k,\sigma_k,\varphi_k}=a_k$, $k=1,2,\ldots,n$. Let $M$ a set of rules, $M=\{{\bf
r}_1,{\bf r}_2, \ldots, {\bf r}_\lambda\}$ then
\[ \left. \begin{array}
{c}
FormalLF[p]\\
\end{array} \right|_{M} =\bigcup_{\nu=1}^\lambda \{ \sum_{\mu =1}^k c_k
L_{1,\mu}^{j_{1,\mu}} \cdots L_{n,\mu}^{j_{n,\mu}} +R \quad \hbox{with} \quad {\cal W}
\to {\bf r}_\nu \in M \}\] The set of substitutions $M$, may be finite or infinite.
\begin{exa}
We shall work with the polynomial $p$ that appears in the example (\ref{ex1}). Let
$M=\{(-2,1,-1)\}$. We determine the following order among the parameters $(W_{3,1,1},
W_{3,2,1}, W_{2,1,2})$. Then $W_{3,1,1}=-2, W_{3,2,1}=1, W_{2,1,2}=-1 $ and $\left.
\begin{array} {c}
FormalLF[p]\\
\end{array} \right|_{M}=11x_1(-2x_1+x_2+x_3)$ $-18x_1(-x_1+x_2)+5x_1+4x_1^2$.
If $M=\{(\varphi,\varphi,\theta)\}, \varphi, \theta \in {\bf R} $ then we set $
W_{3,1,1}=\varphi, W_{3,2,1}=\varphi, W_{2,1,2}=\theta $ and $\left. \begin{array} {c}
FormalLF[p]\\
\end{array} \right|_{M}=11x_1(\varphi x_1+\varphi x_2+x_3)+$
$(-7-11 \varphi)x_1(\theta x_1+x_2)+$ $5x_1+(7\theta -11 \varphi +11 \theta \varphi
)x_1^2$. This is just a re-parameterization of the $FormalLF[p]$.
\end{exa}
In this section we shall examine the positivity of multivariable polynomials via the
Formal-Linear-Like-Factorizations developed above. This is a well-known issue and
there are and other computational approaches - recently, for instance, in
\cite{kn:parriloPHD} - that are mainly based on algorithms that construct sums of squares.
Our contribution essentially relies on the existence of parameters that provide us
with a more flexible tool that is also suitable for solving the positivity problem,
among others. A given polynomial $p$ $\in$ ${\bf R}[x_1,x_2, \ldots, x_n]$, is called
positive, if $p\ge 0$, for every $(x_1,x_2, \ldots, x_n) \in {\bf R}^n$. Let ${\cal
S}$ be the set of values of the parameters $W_{i,\sigma, \varphi}$, such that all but
the square terms of $p$ are eliminated and the coefficients of the remaining terms are
positive real numbers. 
Then, if ${\cal
S}\neq \emptyset$ the polynomial is positive and $\left. \begin{array} {c}
FormalLF[p]\\
\end{array} \right|_{\cal S}$ ia a class of sum of squares of $p$. Here we need to make clear  that all the calculations are
symbolic and not numerical. This permits us to study not only specific polynomials
with numeric coefficients, but also classes of polynomials with parametric
coefficients. In any case, the construction of the set ${\cal S}$ can be carried out
via the following algorithm: \vskip 7 pt \noindent \underline{\bf THE POS-ALGORITHM}
\vskip 7 pt \small {\sf
\par
\noindent {\bf Input:} A multivariable polynomial $p$, a set of undetermined parameters
${\cal W}=\{ W_{i,\sigma,\varphi}\}$, taking values in ${\bf R}$.
\par
\noindent {\bf Initial Conditions:} ${\cal S}=\{\}$, $E=\{ \}$.
\par
\begin{description}
\item{\bf Step 1:} By means of the FORMAL ALGORITHM we get the quantities
$R$, $c_\mu$, $L_{i,\mu}$, $j_{i,\mu}$, $\mu=1, \ldots,k$, $i=1, \ldots, n$. We denote
the coefficients of the remainder $R$ by $c_\mu$, $\mu=k+1, \ldots, k+h$, $h$ the
number of terms of $R$.
\item{\bf Step 2:} {\bf REPEAT FOR } $\mu=1, \ldots, k+h$
\begin{quote}
{\bf IF} some of the exponents $j_{i,\mu}$, $i=1, \ldots, n$ are 
odd numbers {\bf
THEN} $O=O \cup \{ c_{\mu} \}$ {\bf ELSE} $E=E \cup \{ c_{\mu} \}$
\end{quote}
{\bf NEXT $\mu$}
\item{\bf Step 3:} Find the values of the parameters $W_{i,j,\phi}$ that
eliminate the " odd " coefficients and make 
the "even" coefficients positive. In other
words, we construct the set $ S=\{ W_{i,j,\phi}=\lambda_{i,j,\phi} \in {\bf R}: $
$c_{\mu}=0$ and $c_p \ge 0$ with $c_{\mu} \in O$ and $c_p \in E \}$.
\end{description}
 {\bf Output:}
The set ${\cal S}$}. \normalsize
\par
\noindent The proof of the following theorem is straightforward.
\begin{theo}
Let $p \in {\bf R}[x_1, x_2, \ldots, x_n ]$ be a given multivariable polynomial. If
${\cal S}$ is non-void, then $p$ is positive definite, and $\left. \begin{array} {c}
FormalLF[p]\\
\end{array} \right|_{\cal S}$ is a family of \"{} sum of squares \"{} expressions of $p$.
\end{theo}
We would like to make the following remarks in connection with the above algorithm:
\begin{rem}
(i) We can modify the algorithm so that the whole procedure is executed " together "
with the Formal Algorithm, as opposed to after it. This will allow for a quicker
implementation of the method. (ii) Let us suppose that we have a Formal - Linear -
Like - Factorization of a given polynomial $p$. The coefficients of the first terms
with odd exponents contain a small number of parameters ( usually one or two). This
means that they can be eliminated easily for some particular values of the
$W$-parameters. By substituting those values into the other terms we decrease the
number of parameters, thus simplifying the whole computational procedure significantly.
Clearly, the method is not computationally complex and can be carried out normally. Certain algorithms
can be used towards this direction, for instance the Stetter algorithn, to mention but a few.
(iii) If the output of the POS-AlGORITHM is ${\cal S}=\emptyset$, this does not mean
that the polynomial $p$ is not positive or that another sum of squares does not exist.
Nevertheless, our approach can be extended by using nonlinear " factors ", for
instance: $\tilde{L}_{n,k}=W_{n,1,k}x_1+W_{n,2,k}x_2+$ $\cdots +W_{n,(1,2),k}x_1x_2+$
$\cdots +W_{n,(n-1,n-1),k}x_{n-1}^2+x_n^2$. 
In this case we can obtain a further
sums of squares that can successfully deal with cases in which the current
method fails. This will be undertaken
in a future study.
\end{rem}
\begin{exa}
Let us consider the polynomial $p=x^2-2xy+6y^2-4yz+3z^2$. Its
Formal-Linear-Like-Factorization is: $3(z+xW_{3,1,1}+yW_{3,2,1})^2$
$+(y+xW_{2,1,2})(-6W_{3,2,1}-4)$ $(z+xW_{3,1,2}+yW_{3,2,2})$
$+x(6W_{3,2,1}W_{2,1,2}+4W_{2,1,2}-6W_{3,1,1})$ $(z+xW_{3,1,3}+yW_{3,2,3})+$
$(y+xW_{2,1,4})^2$ $(-3W_{3,2,1}^2+6W_{3,2,2}W_{3,2,1}+4W_{3,2,2}+6)+$
$x(y+xW_{2,1,5})$ $(6W_{2,1,4}W_{3,2,1}^2-6W_{3,1,1}W_{3,2,1}+$
$6W_{3,1,2}W_{3,2,1}+6W_{2,1,2}W_{3,2,2}W_{3,2,1}-$ $12W_{2,1,4}W_{3,2,2}W_{3,2,1}$
$-6W_{2,1,2}W_{3,2,3}W_{3,2,1}$ $-12W_{2,1,4}+4W_{3,1,2}+$ $4W_{2,1,2}W_{3,2,2}$
$-8W_{2,1,4}W_{3,2,2}$ $-4W_{2,1,2}W_{3,2,3}$ $+6W_{3,1,1}W_{3,2,3}-2)+R$ (We do not
include the entire remainder because of its size). In order to eliminate the first
non-square term to appear, we set $W_{3,2,1}=-{2 \over 3}$. This transforms the
factorization as follows: $ 3\left( -\frac{2y}{3}+z+xW_{3,1,1} \right)^2$
$-6xW_{3,1,1}(z+xW_{3,1,3}+yW_{3,2,3})$ $+\frac{14}{3}(y+xW_{2,1,4})^2$
$-\frac{2}{3}x(y+xW_{2,1,5})$ $(14W_{2,1,4}-6W_{3,1,1}-9W_{3,1,1}W_{3,2,3}+3)$ $+R'$.
The values $W_{3,1,1}=0$ and $W_{2,1,4}=-\frac{3}{14}$  eliminate the other non-square
terms and the factorization becomes $p=3(-\frac{2y}{3}+z)^2$
$+\frac{14}{3}(-\frac{3x}{14}+y)^2$ $+\frac{11}{14}x^2$. Thus, the set ${\cal S}$ upon
request is ${\cal S}=\{W_{3,2,1}=-\frac{2}{3},W_{3,1,1}=0,W_{2,1,4}=-\frac{3}{14}$
$\quad \hbox{all the other parameters} \quad W_{i,j,\varphi} \quad \hbox{are free} \}$.
The above  " sum of squares " expression of $p$ guarantees that $p$ is positive.
\end{exa}
\section{The Feedback Asymptotic Stabilization}
We are now in a position to apply the entire concept that was raised in earlier
sections to the problem of feedback asymptotically stabilizing a nonlinear system at a
given equilibrium point. This is a well-known topic that has been discussed
extensively in the literature \cite{kn:isidori}, \cite{kn:sontag}. 
In this paper
we adopt a computational approach. Specifically, let us have the continuous nonlinear
system: $\dot{\bf x}={\bf \Phi}({\bf x},{\bf u})$, where ${\bf
x}=(x_1,x_2,\ldots,x_n)^T$ is the state vector, ${\bf u}=(u_1,u_2, \ldots,u_m)^T$ the
input vector and ${\bf \Phi}=[\Phi_1({\bf x},{\bf u}),\ldots,\Phi_n({\bf x},{\bf
u})]^T$, where $\Phi_i({\bf x},{\bf u})$, $i=1,\ldots,n$ are multivariable polynomials
of $(x_1, \ldots, x_n)$ and $(u_1, \ldots, u_m)$ without free terms. Obviously $({\bf
x}^0,{\bf u}^0)=({\bf 0,0})$ is an equilibrium point. Let us denote by $(A,B)$ the
linearization pair of this nonlinear system \"{}around\"{} the origin, (i.e.
$A=\frac{\partial {\bf \Phi}}{\partial {\bf x}}({\bf 0,0}), B=\frac{\partial {\bf
\Phi}}{\partial {\bf u}}({\bf 0,0}))$. As is already known, \cite{kn:sontag}, if the
linear system $(A,B)$is asymptotically controllable (which we can check by testing the
sign of the real parts of the uncontrollable eigenvalues. If all of these are
negative, then the system is asymptotically controllable), then the corresponding
nonlinear system is locally asymptotically stable at the origin. This means that we
can find a matrix $F$ such that the feedback-law constructs a system $\dot{\bf x}={\bf
\Phi}( {\bf x},F({\bf x}))$ which is locally asymptotically stable at the point $({\bf
0,0 })$. The main concern of this paper is to calculate this quantity $F$
computationally, but with the following alterations:
\par
\noindent 1) $F$ does not need to be linear (a matrix) but may also be nonlinear (a
polynomial function). Actually, the problem under examination is that of finding a
state feedback of the form ${\bf u}={\bf a}({\bf x})$, with ${\bf a}({\bf
x})=[a_1({\bf x}),\ldots,a_m({\bf x})]^T$ and $a_i({\bf x})$, $i=1,\ldots,m$
multivariable polynomials too, where the corresponding closed-loop system $\dot{\bf
x}=\bf{ \Phi }({\bf x},{\bf a}({\bf x}))) $ has a global asymptotically stable
equilibrium at $({\bf x}^0,{\bf u}^0)=({\bf 0,0})$.
\par
\noindent 2) Some of the uncontrollable eigenvalues of $(A,B)$ can have zero real
parts.
\par
\noindent 3) The stability is not local but global.
\par
\noindent At this point we introduce the algorithm in order to deal with the feedback
asymptotic stabilization problem. \vskip 10 pt \noindent \small \underline{\bf THE
FEEDBACK-GAS ALGORITHM} \vskip 10 pt \noindent \sf{ {\bf Input:} The polynomials
$\Phi_i({\bf x},{\bf u})$, a polynomial Lyapunov function $L$, the degree of the
feedback law upon request.
\par
\begin{description}
\item{\bf Step 1:} We define the feedback law
${\bf u}={\bf a}({\bf x})$, ${\bf a}({\bf x})=[a_1({\bf x}), a_2({\bf x}), \dots,
a_n({\bf x})]$, with 
\[a_j({\bf x})=\sum_{i_1=1}^nA_{i_1}^{(j)}x_{i_1}
+\sum_{(i_1,i_2)=(1,1)}^{(n,n)}A_{(i_1,i_2)}^{(j)}x_{i_1}x_{i_2}+\]
\[\cdots
+\sum_{(i_1,i_2, \ldots, i_k)=(1,1,\ldots,1)}^{(n,n,\ldots,n)}
A_{(i_1,i_2,\ldots,i_n)}^{(j)} x_{i_1}x_{i_2} \cdots x_{i_n}\]
 and
$A_{(i_1,i_2,\ldots,i_n)}^{(j)}$ unknown parameters taking values in ${\bf R}$.
\item{\bf Step 2:} We define the quantity:
\[V=-\frac{\partial L}{\partial {\bf x}}(\Phi_1 ({\bf x}, 
{\bf a}({\bf x})), \Phi_2
({\bf x}, \ldots, \Phi_n ({\bf x}, {\bf a}({\bf x}))\]
\item{\bf Step 3:} By means of the POS-ALGORITHM we construct the set
${\cal V}$, consisting of those values of the parameters for which $V$ is positive.
\end{description}
{\bf Output:} The set ${\cal V}$. } \rm \normalsize
\begin{theo}
Let us have the nonlinear continuous system: $\dot{\bf x}={\bf \Phi}({\bf x}, {\bf
u})$. Let ${\cal V}$ be the output of the FEEDBACK-GAS-Algorithm. If ${\cal V} \ne
\emptyset$ then the set of feedback laws $\left. \begin{array} {c}
{\bf u}={\bf a}({\bf x})\\
\end{array} \right|_{\cal V}$
make the origin globally asymptotically stable.
\end{theo}
{\bf Proof:} The positive definiteness of the quantity $V$ guarantees that the
Lyapunov function $L$ decreases along the trajectories of the closed-loop system ${\bf
\Phi}({\bf x},{\bf a}({\bf x}))$ and, therefore, the origin is asymptotically stable.
Given that this is the case for every point, the origin is globally asymptotically
stable. $\qed $
\begin{exa}
The angular momentum of a rigid body controlled by two independent torques can
be described, following some simplification, through the equations: $\dot{\bf x}={\bf
\Phi} ({\bf x}, {\bf u})= \left[ \begin{array}{c}
\Phi_1({\bf x}, {\bf u})\\
\Phi_2({\bf x}, {\bf u})\\
\Phi_3({\bf x}, {\bf u})
\end{array} \right]$, with ${\bf x}=[x_1,x_2,x_3]$, ${\bf u}=[u_1(t),u_2(t)]$  and
$\Phi_1=a_1x_2x_3+u_1$, $\Phi_2=a_2x_1x_3+u_2$, $\Phi=a_3x_1x_2$. The quantities
$a_1,a_2$ and $a_3$ are certain constants and $a_3 \ne 0$, (\cite{kn:sontag}, page
176). It can be proved that this system can be globally stabilized about ${\bf x}={\bf
0}$ and ${\bf u}={\bf 0}$. A specific feedback can be constructed by following certain
methods \cite{kn:sontag}.  Following the steps of the Feedback-GAS-Algorithm we select
a Lyapunov function of the form $L=x_1^2+x_2^2+x_3^2$, and a pair of feedback-laws of
the forms: $u_1=A_1x_1+B_1x_2+\Gamma_1x_3+\Delta_1x_2x_3$,
$u_2=A_2x_1+B_2x_2+\Gamma_2x_3+\Delta_2x_1x_3$. Then, we define the quantity:
$V=-\Phi_1(x_1,x_2,x_3,A_1x_1+B_1x_2+\Gamma_1 x_3+\Delta_1x_2x_3)x_1$
$-\Phi_2(x_1,x_2,x_3,A_2x_1+B_2x_2+\Gamma_2 x_3+\Delta_2x_2x_3)x_2$
$-\Phi_3(x_1,x_2,x_3)x_3$. The Formal Algorithm will give the following
Formal-Linear-Like-Factorization of $V$:
\[(a_1+a_2+a_3+\Delta_1+\Delta_2)[W_{2,1,1}
(x_3+x_1W_{3,1,3}+x_2W_{3,2,3})x_1^2\]
\[+(x_2+x_1W_{2,1,7})
(W_{3,1,1}+W_{2,1,1}W_{3,2,1} -2W_{2,1,5}W_{3,2,1}-W_{2,1,1}W_{3,2,3})x_1^2+\]
\[(x_2+x_1W_{2,1,5})^2
(W_{3,2,1}x_1-x_2-x_1W_{2,1,1}) (x_3+x_1W_{3,1,1}+x_2W_{3,2,1})x_1]+\]
\[(\Gamma_2W_{2,1,2}-\Gamma_1)
(x_3+x_1W_{3,1,4}+x_2W_{3,2,4})x_1+ (x_2+x_1W_{2,1,8})\]
\[(-A_2-B_1+2B_2W_{2,1,6}+
\Gamma_2W_{3,1,2}+\Gamma_2W_{2,1,2}W_{3,2,2}-
2\Gamma_2W_{2,1,6}W_{3,2,2}+\Gamma_1W_{3,2,4}-\]
\[\Gamma_2W_{2,1,2}W_{3,2,4})x_1+
(x_2+x_1W_{2,1,6})^2(\Gamma_2W_{3,2,2}-B_2)-\]
\[\Gamma_2(x_2+x_1W_{2,1,2})
(x_3+x_1W_{3,1,2}+x_2W_{3,2,2})+R\] where $R$ is the remainder. The values of the
parameters that eliminate the non-squares terms of the above expressions are
$W_{2,1,1}=0$, $W_{3,2,1}=0$, $W_{3,1,1}=0$, $W_{2,1,4}=\frac{A_2+B_1}{2B_2}$,
$\Gamma_1=0$, $\Gamma_2=0$, $\Delta_1=-\Delta_2-a_1-a_2-a_3$ and
$A_1,A_2,B_1,B_2,\Delta_2$, arbitrary. For those specific values the FormalLF[V]
becomes: $ \frac{(A_2^2+2B_1A_2+B_1^2-4A_1B_2)x_1^2}{4B_2}$ $-B_2\left(
x_2+\frac{A_2+B_1}{2B_2}x_1\right)^2$. In order to get positive coefficients we
further demand $B_2<0$ and $A_1 < \frac{A_2^2+2B_1A_2+B_1^2}{4B_2}$. These conditions
mean that the Lyapunov function decreases along the trajectories of the system, which
guarantees the asymptotic stability of the origin. The family of feedback laws is
given by the relations $u_1=A_1x_1+B_1x_2+\Delta_1x_1x_2$,
$u_2=A_2x_1+B_2x_2+\Delta_2x_1x_3$ with $B_2<0$, $A_1 <
\frac{A_2^2+2B_1A_2+B_1^2}{4B_2}$, $\Delta_1=\Delta_2-a_1-a_2-a_3$ and
$A_2,B_1,\Delta_2$ arbitrary. We can also obtain other classes of feedback, by
choosing other expressions for the quantities $u_1,u_2$ or other Lyapunov
function.

\end{exa}

\begin{exa}Let us suppose that we have the nonlinear system:
$ \dot{\bf x}={\bf \Phi}({\bf x},{\bf u})=\left[\begin{array}{c}
\Phi_1({\bf x},{\bf u})\\
\Phi_2({\bf x},{\bf u}) \end{array} \right] $ with $\Phi_1=4x+8y-{11 \over
16}x^3+5yu_1$$-{5 \over 2}u_1u_2-{4} u_2 $ $+{5}xu_2$, $ \Phi_2=-{ 55 \over
8}x^2-{3}y^3-{1 \over 2}u_2^2x$, where ${\bf x}=[x(t),y(t)]^T$ is the state of the
system, consisting of two functions and ${\bf u}=[u_1(t),u_2(t)]$ the input vector. We
can easily check that the origin is an unstable equilibrium point for the system.
Furthermore, the linearization of this system is not asymptotically controllable,
since the polynomial $\chi_u$ has a root with zero real part (see \cite{kn:sontag}
 for details).
Following the steps of the Feedback-GAS Algorithm, we choose a pair of feedback laws
of the form: $u_1=A_1x+B_1y+\Gamma_1 xy$, $ u_2=A_2x+B_2y+\Gamma_2xy $. Then, we
define the quantity $ V=-\Phi_1(x,y,A_1x+B_1y+\Gamma_1 xy,A_2x+B_2y+\Gamma_2xy )x$
$-\Phi_2(x,y,A_1x+B_1y+\Gamma_1 xy,A_2x+B_2y+\Gamma_2xy)y$. The Formal Algorithm will
give the following Formal-Linear-Like-Factorization of $V$:
\[ FormalLF[V]=-{1 \over 8} \Gamma_2^2W_{2,1,2}^3x^6-{3 \over 2} \Gamma_2^2
W_{2,1,2} \left(y+{1 \over 2} x W_{2,1,2} \right)^2x^4+\]
\[+\frac{-4\Gamma_2A_2^2-20\Gamma_1\Gamma_2A_2-25\Gamma_1^2\Gamma_2}{12B_2}
\left(y+xW_{2,1,10}\right)x^4+{1 \over 2}\Gamma_2^2(y+xW_{2,1,2})^3x^3+\]
\[+{1 \over 32}(-B_2^4+32A_2B_2+80\Gamma_1B_2-160\Gamma_1+80\Gamma_2B_1)\cdot\]
\[\cdot \left(y-x\frac{B_2^6-576A_2^2-2880A_2\Gamma_1-2880A_1\Gamma_2+5760\Gamma_2}
{72(-B_2^4+32A_2B_2+80\Gamma_1B_2-160\Gamma_1+80\Gamma_2B_1)}\right)^2x^2+\]
\[+\Gamma_2 \left( y -x\frac{(-2A_2-5\Gamma_1)}{6B_2} \right)^3B_2x^2 +{5 \over 2}(B_1B_2-2B_1)\cdot\]
\[ \cdot \left(y-x\frac{(-20B_2A_1+40A_1-32\Gamma_2-20A_2B_1+40B_2-55)}{40B_1(B_2-2)}\right)^2x+\]
\[+4(B_2-2)(y+xW_{2,1,13})x+3\left(x{B_2^2 \over 24}+y\right)^4+R\]
where $R$ is the remainder (not included because of its size). The values of the
parameters that eliminate the non square terms of the above expression and make the
coefficients of the remaining terms positive, are ${\cal V}=\{ A_1=2,$
$A_2=\frac{5}{4B_1},$ $B_2=2,$ $\Gamma_2=0,$ with $0<B_1<{5 \over 4}$ and
$\frac{16B_1^3-225B_1-32\sqrt{111}\sqrt{5B_1^5-B_1^6}} {900B_1^2} <\Gamma_1 <$ $
\frac{16B_1^3-225B_1+32\sqrt{111}\sqrt{5B_1^5-B_1^6}} {900B_1^2}\}$ and
$W_{i,j,\sigma}$ arbitrary. By taking, for instance, $B_1=1, \Gamma_1={1\over 2}$ the
$FormalLF[V]$ becomes: 
\[\left. \begin{array} {c}
FormalLF[V]\\
\end{array} \right|_{\cal V}=
\frac{20375}{663552}x^4 +2\left(\frac{659}{1152}x+y\right)^2x^2+
x^2+3\left(\frac{x}{6}+y\right)^4\]
 which means that the Lyapunov function
decreases along the trajectories of the system. This guarantees the asymptotic
stability of the origin. The family of the feedback laws is given by the relations
$u_1=x+B_1y+\Gamma_1xy$, $u_2=\frac{5}{4B_1}x+2y$, with $0 < B_1< {5 \over 4}$, and
$\frac{16B_1^3-225B_1-32\sqrt{111}\sqrt{5B_1^5-B_1^6}} {900B_1^2} <\Gamma_1 <$ $
\frac{16B_1^3-225B_1+32\sqrt{111}\sqrt{5B_1^5-B_1^6}} {900B_1^2}$. It should be noted
that the values $0$ and ${5 \over 4}$ for the parameter $B_1$ are a kind of
bifurcation point for the feedback expression, above or below which the asymptotic
stabilization is not achieved.

\end{exa}


\begin{thebibliography}{99}


\bibitem{kn:lasserre} J.B.Lasserre " Global optimization with polynomials
 and the problem of moments. " SIAM Journal on Optimization, Vol. 11, N0. 3, 
 pp. 796-817, 2001.
 

\bibitem{kn:choi} M.D.Choi and T. Y. Lam, and B. Reznick " Sums of squares of real polynomials. "
Proceedings of Symposia in Pure Mathematics, 58 (2): 103-126, 1995.

\bibitem{kn:hasan} M.A.Hasan and A.A.Hasan " A procedure for the positive definiteness
of forms of even order. " IEEE Transactions on Automatic Control, 41(4): 615-617, 1996.


\bibitem{kn:hochbaum} D.S.Hochbaum, editor. " Approximation algorithms for NP-hard problems. "
PWS Publishing Company, 1997.

\bibitem{kn:isidori} A. Isidori. " Nonlinear Control Systems. " Springer - Verlag, Berlin, third
edition, 1997.

\bibitem{kn:mishra} B.Mishra. " Algorithmic Algebra. " Springer-Verlag, 1993.

\bibitem{kn:powers} V.Powers and T.W\"{o}rmann. " An algorithm for sums of squares of real polynomials. "
Journal of Pure and Applied Algebra,127, 99-104,1998.

\bibitem{kn:boyd} L.Vandenberghe and S.Boyd. " Semidefinite Programming ". SIAM Review, 38(1), 49-95,
March 1996.


\bibitem{kn:cox} D. Cox, J. Little, D. O'Shea. (1997). "Ideals, Varieties
and Algorithms". {\it Springer-Verlag, New York}.


\bibitem{kn:parriloPHD} P.A.Parrilo (2000) " Structured semidefinite programs
and semialgebraic geometry methods in robustness and optimization ". PhD thesis,
California Institute of Technology.

\bibitem{kn:parrilo} P.A.Parrilo (2002) " An explicit construction of distinguished representations
of polynomials nonnegative over finite sets ". IFA Technical Report AUT02-02.



\bibitem{kn:prestel} A.Prestel and C.N.Delzell, (2001) " Positive polynomials:
from Hilbert's 17th problem to real algebra. " Springer Monographs in Mathematics.
Springer, New York.



\bibitem{kn:karkanias} N.Munro (editor), (1999) " Symbolic Methods in control System Analysis and Design. "
IEE Control Engineering Series, 56.

\bibitem{kn:sontag} E.Sontag (1997) " Mathematical Control Theory ", Texts in Applied Mathematics,
Springer, New York.


\bibitem{kn:hercma} S.Kotsios (2003) " The problem 
of positive defineteness through a formal factorization of polynomials. "
 HERCMA - Congress 2003, ATHENS.

\end{thebibliography}
\end{document}